\documentstyle [11pt,epsf, graphicx, amssymb]{article}

\begin{document}

\def\l{\lambda}
\def\m{\mu}
\def\a{\alpha}
\def\b{\beta}
\def\g{\gamma}
\def\d{\delta}
\def\e{\epsilon}
\def\o{\omega}
\def\O{\Omega}
\def\v{\varphi}
\def\t{\theta}
\def\r{\rho}
\def\bs{$\blacksquare$}
\def\bp{\begin{proposition}}
\def\ep{\end{proposition}}
\def\bt{\begin{th}}
\def\et{\end{th}}
\def\be{\begin{equation}}
\def\ee{\end{equation}}
\def\bl{\begin{lemma}}
\def\el{\end{lemma}}
\def\bc{\begin{corollary}}
\def\ec{\end{corollary}}
\def\pr{\noindent{\bf Proof: }}
\def\note{\noindent{\bf Note. }}
\def\bd{\begin{definition}}
\def\ed{\end{definition}}
\def\C{{\mathbb C}}
\def\P{{\mathbb P}}
\def\Z{{\mathbb Z}}
\def\d{{\rm d}}
\def\deg{{\rm deg\,}}
\def\deg{{\rm deg\,}}
\def\arg{{\rm arg\,}}
\def\min{{\rm min\,}}
\def\max{{\rm max\,}}

\newcommand{\norm}[1]{\left\Vert#1\right\Vert}
\newcommand{\abs}[1]{\left\vert#1\right\vert}

\newcommand{\set}[1]{\left\{#1\right\}}
\newcommand{\setb}[2]{ \left\{#1 \ \Big| \ #2 \right\} }

\newcommand{\IP}[1]{\left<#1\right>}
\newcommand{\Bracket}[1]{\left[#1\right]}
\newcommand{\Soger}[1]{\left(#1\right)}

\newcommand{\Integer}{\mathbb{Z}}
\newcommand{\Rational}{\mathbb{Q}}
\newcommand{\Real}{\mathbb{R}}
\newcommand{\Complex}{\mathbb{C}}

\newcommand{\eps}{\varepsilon}
\newcommand{\To}{\longrightarrow}
\newcommand{\varchi}{\raisebox{2pt}{$\chi$}}

\newcommand{\E}{\mathbf{E}}
\newcommand{\Var}{\mathrm{var}}

\def\squareforqed{\hbox{\rlap{$\sqcap$}$\sqcup$}}
\def\qed{\ifmmode\squareforqed\else{\unskip\nobreak\hfil
\penalty50\hskip1em\null\nobreak\hfil\squareforqed
\parfillskip=0pt\finalhyphendemerits=0\endgraf}\fi}

\renewcommand{\th}{^{\mathrm{th}}}
\newcommand{\Dif}{\mathrm{D_{if}}}
\newcommand{\Difp}{\mathrm{D^p_{if}}}
\newcommand{\GHF}{\mathrm{G_{HF}}}
\newcommand{\GHFP}{\mathrm{G^p_{HF}}}
\newcommand{\f}{\mathrm{f}}
\newcommand{\fgh}{\mathrm{f_{gh}}}
\newcommand{\T}{\mathrm{T}}
\newcommand{\K}{^\mathrm{K}}
\newcommand{\PghK}{\mathrm{P^K_{f_{gh}}}}
\newcommand{\Dig}{\mathrm{D_{ig}}}
\newcommand{\for}{\mathrm{for}}
\newcommand{\End}{\mathrm{end}}

\newtheorem{th}{Theorem}[section]
\newtheorem{lemma}{Lemma}[section]
\newtheorem{definition}{Definition}[section]
\newtheorem{corollary}{Corollary}[section]
\newtheorem{proposition}{Proposition}[section]

\begin{titlepage}

\begin{center}

\topskip 5mm

{\LARGE{\bf {Norming Sets and Related Remez-type Inequalities}}}

 \vskip 5mm

{\large {\bf  A. Brudnyi $^{*}$, Y. Yomdin$^{**}$}}

\end{center}

{$^{*}$ Department of Mathematics and Statistics, University of
Calgary, Calgary, Alberta, T2N 1N4 Canada. E-mail:
abrudnyi@ucalgary.ca}

{$^{**}$ Department of Mathematics, The Weizmann Institute of
Science, Rehovot 76100, Israel. E-mail:
yosef.yomdin@weizmann.ac.il}

\begin{center}

{ \bf Abstract}
\end{center}

{\small The classical Remez inequality (\cite{Rem}) bounds the maximum of the 
absolute value of a real polynomial $P$ of degree $d$ on $[-1,1]$
through the maximum of its absolute value on any subset $Z\subset [-1,1]$ of
positive Lebesgue measure. Extensions to several variables and
to certain sets of Lebesgue measure zero, massive in a much weaker
sense, are available (see, e.g., \cite{Bru.Gan,Yom5,Bru}).

Still, given a subset $Z\subset [-1,1]^n\subset {\mathbb R}^n$ it is not easy to determine
whether it is ${\mathcal P}_d({\mathbb R}^n)$-norming (here ${\mathcal P}_d({\mathbb R}^n)$ is the space of real polynomials of degree at most $d$ on ${\mathbb R}^n$), i.e. satisfies a Remez-type inequality:
$\sup_{[-1,1]^n}|P|\le C\sup_{Z}|P|$ for all $P\in {\mathcal P}_d({\mathbb R}^n)$ with $C$ independent of $P$. (Although
${\mathcal P}_d({\mathbb R}^n)$-norming sets are exactly those not contained in any algebraic hypersurface
of degree $d$ in ${\mathbb R}^n$, there are many apparently unrelated
reasons for $Z \subset [-1,1]^n$ to have this property.)

In the present paper we study norming sets and related Remez-type inequalities in a general
setting of finite-dimensional linear spaces $V$ of
continuous functions on $[-1,1]^n$, remaining in most of the examples in
the classical framework. First, we discuss some sufficient conditions for $Z$ to be
$V$-norming, partly known, partly new, restricting ourselves to the simplest 
non-trivial examples. Next, we extend the Turan-Nazarov inequality for exponential 
polynomials to several variables, and on this base prove a new fewnomial Remez-type inequality. 
Finally, we study the family of optimal constants $N_{V}(Z)$ in the Remez-type inequalities for $V$, 
as the function of the set $Z$, showing that it is Lipschitz in the Hausdorff metric.}

\medskip

\smallskip


\begin{center}
------------------------------------------------
\end{center}
Research of the first author was partially supported by NSERC.\newline
Research of the second author was supported by the ISF, Grant No. 779/13 and by
the Minerva foundation.

\end{titlepage}

\newpage

\section{Introduction}
\setcounter{equation}{0}

The classical Chebyshev inequality, see, e.g., \cite[pp.~67--68]{Tim}, bounds the maximum of
the absolute value of a polynomial $P$ of degree $d$ on
$[-1,1]$ through the maximum of its absolute value on an interval $[a,b]\subset [-1,1]$.
In a little known and hardly available paper \cite{Rem}, Remez generalized the Chebyshev inequality by
replacing $[a,b]$ by an arbitrary measurable subset $Z\subset [-1,1]$:

\medskip

\bt\label{Rem} Let $P\in {\mathbb R} [x]$ be a polynomial of degree $d$. Then for
any measurable $Z\subset [-1,1]$
\be\label{Rem.in}
\sup_{[-1,1]} \vert P\vert \leq T_d\left({{4-\mu}\over {\mu}}\right)\sup_Z \vert P \vert ,
\ee
where $\mu=\mu_1(Z)$ is the Lebesgue measure of $Z$ and
$T_d$ is the Chebyshev polynomial of degree $d$.
\et

\medskip

This result has been rediscovered several times, see, e.g.,
\cite{DR} and \cite[Lm.~2]{Bru.Gan}, but the Remez
proof is still the most simple and elegant.

\medskip

A multidimensional inequality of this kind is due to Yu. Brudnyi and Ganzburg \cite{Bru.Gan}:

\bt\label{Rem.mult} Let $
B\subset {\mathbb R}^n$ be a convex body and $Z\subset
B$ be a measurable subset. Then for every real polynomial
$P$ on ${\mathbb R}^n$ of degree $d$,
\be\label{BG} \sup_{
B}\vert P \vert \leq T_d \left({{1+(1-\lambda)^{1\over n}}\over
{1-(1-\lambda)^{1\over n}}}\right) \sup_{Z}\vert P \vert.\ee Here
$\lambda= {{\mu_n(Z)}\over {\mu_n(B)}}$ with $\mu_n$ being
the Lebesgue measure on ${\mathbb R}^n$.
\et

\medskip

This inequality
is sharp and for $n=1$ coincides with (\ref{Rem.in}).

\medskip

It is well known that inequalities of the form (1.2)
may be true also for some sets $Z$ of Lebesgue measure zero and even for
certain finite sets $Z$, see, e.g.,
\cite{Ber2,Bru,Bru.Bru, Cop.Riv,Fav, Fri.Yom, Rah,Yom5,Yom6,Yom7,Zer1}.

\medskip

There are numerous generalizations of the above inequalities (referred to as Remez-type inequalities) to wider classes of
functions. Recently there has been a considerable interest
in such inequalities in connection with various problems of analysis, see, e.g., the Introduction to \cite{Bru2} and references therein,
results and references in \cite[Ch.~2,~9,~10]{BB},  \cite{And, BE, BDSV, Erd, Kro}, etc. Some of the results below
can be extended in an appropriate form to these classes. 

\medskip  In the present paper we study the problem of characterizing the objects subject to the following definition.

Let $V\subset C(Q^n_1)$ be a finite-dimensional subspace of real continuous functions on the closed unit cube $Q^n_1\subset {\mathbb R}^n$. 
\bd\label{definite}
A compact subset $Z\subset Q^n_1$ is said to be $V$-norming if there exists a constant $C>0$ such that for every $f\in V$
\be\label{FRemez}
\max_{Q_1^n}|f|\le C\cdot\max_{Z} |f|.
\ee
The minimum of all such constants $C$ is denoted by $N_{V}(Z)$ and is called the $V$-norming constant of $Z$.
\ed
(The notion originates from the Banach Space theory: the family $\{\delta_z\}_{z\in Z}\subset V^*$ of evaluation functionals at points of $Z$ is a norming set for $V$, i.e.
$\|f\|_Z:=\sup_{z\in Z}|\delta_z(f)|$, $f\in V$, is a norm equivalent to the supremum norm on $V$.)

\medskip

Even if $V:={\mathcal P}_d({\mathbb R}^n)$, the space of real polynomials of degree at most $d$ on ${\mathbb R}^n$, the general problem of characterizing
sets $Z$ for which Remez-type inequality (\ref{FRemez}) is valid, remains generally open. As we will see below, there is a wide variety of apparently unrelated geometric, algebraic, arithmetic, etc. {\it sufficient conditions} on such $Z$, which are difficult to present in
a coherent way.


In principle, there is a very simple description of $V$-norming sets:

\bp\label{prop1} A compact subset $Z\subset Q^n_1$ is $V$-norming
(equivalently, $N_{V}(Z)< \infty$) if and
only if the space $V|_Z$ of restrictions of functions in $V$ to $Z$ is of dimension ${\rm dim}\, V$. This is, in turn, equivalent to the condition that $Z$ is not contained in any zero-level set $\{x\in Q^n_1\, :\, f(x)=0\}, \ f\in V$.
\ep
\pr Indeed, (\ref{FRemez}) implies that $Z$ is $V$-norming if and only if the restriction map $r: V\rightarrow V|_Z$, $r(f):=f|_{Z}$,
is an isomorphism of Banach spaces equipped with the corresponding supremum norms.\qquad
$\square$

\smallskip

However, in general it is not easy to reformulate this
condition in an ``effective'' way and to provide explicit
bounds on $N_V(Z)$ starting with an explicitly given $Z$.

\smallskip

In the present paper we first discuss (in Section \ref{Sec2}) some sufficient conditions for $Z$ to be ${\mathcal P}_d({\mathbb R}^n)$-norming, partly known, partly new. This includes ``massiveness'', algebraic, and topological properties of $Z$. Next, in Section \ref{Few} we extend the Turan-Nazarov inequality for exponential 
polynomials to several variables, and on this base prove a  new fewnomial Remez-type inequality. 
Finally, in Section \ref{Sec3} we study the behavior of the best constant $N_V(Z)$ in the Remez-type inequality, as the function of the set $Z$, showing, in particular, that it is Lipschitz in the Hausdorff metric. 

\medskip

\noindent {\em Acknowledgment.} We thank Len Bos and Yuri Brudnyi for useful comments improving the presentation of the paper.

\section{Examples of Remez Sets}\label{Sec2}
\setcounter{equation}{0}

As it was mentioned in the Introduction, the problem of characterizing
(in geometric terms) those sets $Z$ for which Remez-type inequality is valid is  generally open.
In this section we provide a small number of a wide variety of apparently unrelated geometric, algebraic, arithmetic, etc. {\it sufficient conditions} for $Z$ to be $V$-norming. We restrict ourselves to the simplest non-trivial examples of each kind, mostly for families $V$ of real polynomials or analytic functions.

\subsection{Interpolation System}\label{Intrpol}

For finite
sets $Z$ it is possible (in principle) to write an explicit answer
through the determinants arising in the corresponding interpolation systems.

\smallskip

Indeed, suppose ${\mathcal F}:=\{f_1,\dots, f_l\}\subset V$, $l:={\rm dim}\, V$, is a basis and our set
$Z$ contains exactly $l$ points $x^1,\dots,x^l \in Q^n_1$.
Assuming that the values of $f=\sum_{i=1}^l a_i f_i\in V$ on $Z$ are given, $f(x^j)=v_j, \ j=1,\dots,l,$
we get the following interpolation system:

\be \sum_{i=1}^l a_i f_i(x^j)= v_j,
\qquad j=1,\dots,l.
\ee
Considered as a linear system with respect to the unknown variables
$a_i$, this is a multi-dimensional Vandermonde-like system with the matrix
$M_Z=\bigl(f_i(x^j)\bigr)$. It is uniquely solvable if and only if its determinant
$\Delta_{\mathcal F}(x^1,\dots,x^l)= \det M_Z$ is nonzero. The $V$-norming constant $N_V(Z)$ is
exactly the norm of the inverse matrix $M_Z^{-1},$ considered as the operator
from the space of functions on $Z$ to the space
$V$ both equipped with the corresponding supremum norms. An easy
application of Cramer's rule gives us the bound of the form
\be\label{cramer}
N_V(Z)\leq {\bigl({\max_{1\le i\le l}\left\{\max_{Q^n_1}|f_i|\right\}\bigr)^l\cdot l\cdot l!}\over {\left|\Delta_{\mathcal F}(x^1,\dots,x^l)\right|}}.
\ee

Thus we have:

\bp
A set $Z=\{x^1,\dots,x^l\}\subset Q^n_1$ is $V$-norming
if and only if $\Delta_{\mathcal F}(x^1,\dots,x^l)\ne 0$. In this case, the upper bound for
$N_V(Z)$ is given by (\ref{cramer}).
\ep
For specific families ${\mathcal F}\subset V$ and sets $Z$ more accurate
estimates of $N_V(Z)$ in terms of the characteristics of the matrix $M$ can be produced.

\medskip

In this extremal setting the Remez-type inequality is essentially equivalent to
the stability estimate of the multidimensional interpolation problem (see, e.g.,
\cite{Lor} and references therein). As an immediate consequence
we get:

\bc
For each compact $V$-norming set $Z\subset Q^n_1 ,$ the norming constant $N_{V}(Z)$ satisfies
\be\label{bound1}
N_{V}(Z)\leq \inf_{Z'\subset Z,\, \#Z'=l} N_{V}(Z')\le \inf_{{\mathcal F}\subset V}{\bigl({\max_{1\le i\le l}\left\{\max_{Q^n_1}|f_i|\right\}\bigr)^l\cdot l\cdot l!}\over {\sup_{x^1,\dots, x^l\in Z}\left|\Delta_{\mathcal F}(x^1,\dots,x^l)\right|}}.
\ee
(Here $\# Z'$ stands for the cardinality of the subset $Z'$ and $\mathcal F$ runs over all bases in $V$.)
\ec

The first inequality in (\ref{bound1}) is almost optimal as the following result shows.

\bp
For each compact $V$-norming set $Z\subset Q^n_1 ,$ there exists a subset $Z'\subset Z$ of cardinality $l$ such that
\[
\frac{1}{l}\cdot N_{V}(Z')\le N_{V}(Z).
\]
\ep
Thus,
\be\label{equiv}
\frac{1}{l}\cdot \inf_{Z'\subset Z,\, \#Z'=l} N_{V}(Z')\le N_{V}(Z)\leq \inf_{Z'\subset Z,\, \#Z'=l} N_{V}(Z').
\ee

\pr For a fixed basis $\mathcal F$ in $V$ let $Z'=\{x_*^1,\dots x_*^l\}\subset Z$ be such that
\[
\left|\Delta_{\mathcal F}(x_*^1,\dots,x_*^l)\right|=\sup_{x^1,\dots, x^l\in Z}\left|\Delta_{\mathcal F}(x^1,\dots,x^l)\right|.
\]
(Such points exist because $Z$ is compact and $\Delta_{\mathcal F}$ is a continuous function on $(Q^n_1)^l\subset {\mathbb R}^{nl}$.)
Since ${\rm dim}\, V|_{Z}=l$, and evaluations $\delta_z$ at points $z\in Z$ determine bounded linear functionals on $V|_{Z}$, the Hahn-Banach theorem implies easily that ${\rm span}\,\{\delta_z\}_{z\in Z}= (V|_{Z})^*$ and, hence, $\Delta_{\mathcal F}(x_*^1,\dots,x_*^l)\ne 0$.
Next, we define functions $L_i\in V_{\mathcal F}$ by the formulas
\be\label{eq2}
L_i(x):=\frac{\Delta_{\mathcal F}(x_*^1,\dots, x_*^{i-1}, x, x_*^{i+1},\dots, x_*^l)}{\Delta_{\mathcal F}(x_*^1,\dots, x_*^l)},\quad x\in Q^n_1,\quad 1\le i\le n.
\ee
Clearly, they satisfy the following properties
\be\label{eq2'}
L_i(x_*^j)=\delta_{ij}\quad {\rm (the\ Kronecker-delta)}\quad {\rm and}\quad \max_Z |L_i|\le 1.
\ee
For a function $h$ defined on $Z'$ the Lagrange interpolation is given by the formula
\be\label{eq3}
(Lh)(x):=\sum_{i=1}^l h(x_*^i)L_i(x),\quad x\in Q^n_1.
\ee
From (\ref{eq2'}) and (\ref{eq3}) we obtain
\[
N_{V}(Z')\le \sum_{i=1}^l \max_{Q^n_1}|L_i|\le l\cdot N_{V}(Z),
\]
as required.\qquad $\Box$

\subsection{Sets with Algebraically Independent Coordinates}

Assume that $Z=\{x^1,\dots,x^s\}\subset Q^n_1\subset {\mathbb R}^n$, where $s=\left(^n_d\right)\, (={\rm dim}\, {\mathcal P}_d({\mathbb R}^n))$.
As it was shown above, $Z$ is ${\mathcal P}_d({\mathbb R}^n)$-norming (i.e., satisfies Remez-type inequality (\ref{FRemez}) for real polynomials on ${\mathbb R}^n$ of degree at most $d$)  if and only if the Vandermonde matrix $M_Z$
determined with respect to the basis ${\mathcal M}_{d,n}$ of monomials in ${\mathcal P}_d({\mathbb R}^n)$ is nondegenerate, i.e. its determinant
$\Delta_{{\mathcal M}_{d,n}}(x^1,\dots,x^s)\neq 0$. But $\Delta_{{\mathcal M}_{d,n}}$ is a polynomial with integer coefficients in the coordinates of $x^1,\dots,x^s$. Therefore, if they are algebraically independent over $\mathbb Q$, then $\Delta_{{\mathcal M}_{d,n}}(x^1,\dots,x^s)\neq 0$. For instance, due to the classical Lindemann-Weierstrass theorem \cite{W}, the latter is true if all these coordinates are exponents of linearly independent over ${\mathbb Q}$ algebraic numbers.
Presumably, in some specific examples (e.g., if the coordinates of $x^1,\dots,x^s$ are Liouville numbers) the norming constant $N_{{\mathcal P}_d({\mathbb R}^n)}(Z)$ can be estimated explicitly.

\subsection{Hausdorff Measure and Metric Entropy}

By Theorem 1.2 above each compact subset $Z\subset Q^n_1$ of positive Lebesgue
$n$-measure is ${\mathcal P}_d({\mathbb R}^n)$-norming for each natural $d$ and
\be\label{meas}
N_{{\mathcal P}_d({\mathbb R}^n)}(Z)\leq T_d \left({{1+(1-\lambda)^{1\over n}}\over
{1-(1-\lambda)^{1\over n}}}\right)<\left(\frac{4n}{\lambda}\right)^d
\ee
with $\lambda = \mu_n(Z)$. Similarly, if $V$ consists of real analytic functions defined in a neighbourhood of $Q^n_1$ and $Z$ is as above, then
it is $V$-norming and
\[
N_{V}(Z)\le E \left({{1+(1-\lambda)^{1\over n}}\over
{1-(1-\lambda)^{1\over n}}}\right)^C< \left(\frac{4n}{\lambda}\right)^C,
\]
where $E(x):=x+\sqrt{x^2-1}$, $|x|\ge 1$, and $C$ is a constant depending on $V$ only.

This follows from an inequality similar to (\ref{BG}) for real analytic functions, see \cite{Bru-1},  \cite{Bru0}.

Next, it is shown in \cite{Bru.Bru} that compact subsets $Z\subset Q^n_1$ of Hausdorff
dimension greater than $n-1$ are ${\mathcal P}_d({\mathbb R}^n)$-norming for each $d$.

\medskip

In \cite{Yom5,Bru} some ${\mathcal P}_d({\mathbb R}^n)$-norming sets have been characterized in terms
of their metric entropy. Let us recall that the covering number $M(\e,X)$
of a compact metric space $X$ is the minimal number of closed $\e$-balls covering $X$
(see \cite{Kol.Tih}). Below $M(\e,X)$ are defined for compact subsets $X\subset {\mathbb R}^n$ equipped with the
induced $l^\infty$ metric, that is closed $\e$-balls in this metric are intersections with $X$ of closed cubes of sidelength $2\e$ with centers at points of $X$.

\smallskip

\bd \label{span}
Let $Z$ be a compact subset of $Q^n_1$. The
metric $(d,n)$-span $\omega_d(Z)$ of $Z$ is defined
as
\be
\omega_{d,n}(Z): = \sup_{\e>0} {\e}^n[M(\e,Z)- M_{n,d}(\e)].
\ee
\ed
Here $M_{n,d}(\e):=\sum_{i=0}^{n-1}C_i(n,d)\left({1\over \e}\right)^i$ is a universal polynomial of degree $n-1$ in $\frac{1}{\e}$ whose coefficients are positive
numbers related to Vitushkin's bounds for covering numbers of polynomial sub-level sets of degree $d$ (see \cite{Vit1,Yom5,Fri.Yom1}). The explicit formula for
$M_{n,d}$ is given in \cite{Yom5}. In particular,
\[
M_{1,d}(\e)=d, \qquad M_{2,d}(\e)=(2d-1)^2 + 8d\cdot\left({1\over \e}\right).
\]

The following result has been established in \cite{Yom5}:

\bt\label{Discr.Rem}
If $\omega_{d,n}(Z) = \o >0$, then $N_{{\mathcal P}_d({\mathbb R}^n)}(Z)<\infty$ and satisfies
\be
N_{{\mathcal P}_d({\mathbb R}^n)}(Z)\leq T_d \left({{1+(1-\o)^{1\over n}}\over
{1-(1-\o)^{1\over n}}}\right)\quad (=: R_d(\o)).
\ee
\et
Thus, in some cases the Lebesgue measure $\mu_n(Z)$
in Theorem 1.2 can be replaced with $\omega_{d,n}(Z)$.

Note that the metric $(d,n)$-span
$\omega_{d,n}(Z)$ may be positive even for some finite sets. Consider,
for example, finite subsets of ${\mathbb R}$.

\bc \label{Discr.Rem1} A set $Z=\{x_1,\dots,x_m\} \subset [-1,1], \ x_i\ne x_j$ for all $i\ne j$, is
${\mathcal P}_{d}({\mathbb R})$-norming if and only if $m\geq d+1$. In this
case $\omega_{d,1}(Z)\geq \delta$ where $\delta$ is the minimal distance
between distinct $x_i$ and $x_j$ in $Z$, and $N_{{\mathcal P}_{d}({\mathbb R})}(Z)\leq R_d(\delta)=
T_d \left({{2-\delta}\over {\delta}}\right).$\ec\pr It is enough to take $\e>0$
tending to $\delta$ from the left in the definition of $\omega_{d,1}(Z)$, and use Theorem 2.1.\qquad
$\square$

\medskip

Some other examples are given in \cite{Yom5}. In particular, $\omega_{d,n}(Z)>0$
for each $Z$ with the entropy (or box) dimension $\dim_e Z$ greater than $n-1$.

For similar results
for spaces $V$ of real analytic functions, see \cite{Bru}.

\smallskip

Still, subsets $Z \subset Q^n_1$ with $\omega_{d,n}(Z)>0$ are
(in a certain discrete sense) ``massive in dimension $n-1$". Below
we give examples of ${\mathcal P}_d({\mathbb R}^n)$-norming sets $Z$ in $Q^n_1$ which are
contained in certain analytic curves in ${\mathbb R}^n$, and which have
$\omega_{d,n}(Z)=0$. So ``massiveness" is just one of possible
geometric reasons for a set to be ${\mathcal P}_d({\mathbb R}^n)$-norming.

\subsection{Capacity}
Another class of ``massive'' norming sets consists of the, so-called, nonpluripolar subsets of ${\mathbb R}^n$.
Recall that a compact subset $Z\subset Q^n_1$ is {\em pluripolar} if there exists a non-identically $-\infty$ plurisubharmonic
function $u$ on ${\mathbb C}^n$ such that $u|_Z\equiv-\infty$. It is known (see, e.g., \cite{Kli}) that a compact subset $Z\subset Q^n_1$ is nonpluripolar if and only if
there exists a constant $C>0$ depending on $Z$ and $n$ only such that $Z$ is ${\mathcal P}_d({\mathbb R}^n)$-norming for all $d$ and the norming constants satisfy
\[
N_{{\mathcal P}_d({\mathbb R}^n)}(Z)\leq C^d.
\]
For instance, inequality (\ref{meas}) shows that any compact subset $Z\subset Q^n_1$ of positive Lebesgue $n$-measure is  nonpluripolar. However, nonpluripolar sets in ${\mathbb R}^n$
may be of an arbitrary small positive Hausdorff dimension (e.g., the $n$-fold direct product of Cantor sets in $[-1,1]$ of sufficiently small positive Hausdorff dimensions is nonluripolar in
${\mathbb R}^n$, see also \cite{Lab}.) If $Z\subset Q^n_1$ is nonpluripolar, then the upper bound for $N_{{\mathcal P}_d({\mathbb R}^n)}(Z)$ can be expressed also in terms of capacity ${\rm cap}(Z)$ of $Z$, a positive number defined in one of the
following equivalent ways: in terms of the Monge-Amp\`{e}re measure of $Z$, the Robin constant of $Z$, the Chebyshev constant of $Z$ or the transfinite diameter of $Z$, see \cite{AT}, \cite{Kli}, \cite{TL}. Then for such $Z$ and all $d\in {\mathbb N}$ one has
\[
\ln\left(N_{{\mathcal P}_d({\mathbb R}^n)}(Z)\right)\le \frac{c\cdot d}{{\rm cap}(Z)},
\]
where $c>0$ depend on $n$ only, and ${\rm cap}(Z)$ is defined in terms of the Monge-Amp\`{e}re measure of $Z$, see \cite{AT}.

Finally, observe that due to Proposition 1.1 each nonpluripolar compact subset $Z\subset Q^n_1$ is $V$-norming for every finite-dimensional space of real analytic functions defined in a neighbourhood of $Q^n_1$.

\subsection{Nodal Sets of Elliptic PDEs}
We consider a homogeneous elliptic differential equation of the form
\be\label{ellipt}
{\mathcal L}u\equiv\sum_{i,j=1}^n a_{ij}(x)\partial_{ij}u+\sum_{i=1}^n b_i(x)\partial_i u+ c(x)u=0
\ee
defined in an open neighbourhood of $Q^n_1$, where the coefficients $a_{ij}$ satisfy
\[
\sum_{i,j=1}^n a_{ij}(x)\xi_i \xi_j\ge \lambda \|\xi\|_{2}^2\quad {\rm for\ any}\quad \xi\in {\mathbb R}^n,\ x\in Q^n_1
\]
for some positive constant $\lambda$. We assume that $a_{ij}$ are Lipschitz and $b_i$ and $c$ are at least
bounded. The Lipschitz condition for the leading coefficients is essential. It implies
the unique continuation for the operator $\mathcal L$. In other words, if a solution $u$ vanishes
to an infinite order at a point in $Q^n_1$, then $u$ is identically zero, see \cite{Ar}.

For any $C^2$ non-identically zero solution $u$ in $Q^n_1$ we define the nodal set
\[
{\mathcal N}(u):=\{x\in Q^n_1\, :\, u(x)=0\}.
\]
According to \cite{HS}, the set ${\mathcal N}(u)$ has finite $(n-1)$-dimensional Hausdorff measure. Thus from here and Proposition 1.1 we obtain:

{\em Suppose that $V$ is a finite-dimensional space of $C^2$ solutions of equation (\ref{ellipt}) and $Z\subset Q^n_1$ is compact of infinite $(n-1)$-dimensional Hausdorff measure. Then $Z$ is $V$-norming.}

\subsection{Algebraic Curves of High Degree}

One apparent algebraic-geometric reason for a set to be ${\mathcal P}_d({\mathbb R}^n)$-norming
is that algebraic sets $Z$ of degree higher than $d$ ``generically" cannot
be contained in a hypersurface of degree $d$. There are plenty of ways
in which this general claim can be transformed into a Remez-type inequality.
We give here only one simple example, where computations are fairly
straightforward.

Consider a curve $S \subset Q^n_1 \subset {\mathbb R}^n$ given in
parametric form by $x=\Psi(t)$ where $\Psi$ is defined by

\be x_1=t^{d_1}, \ x_2=t^{d_2},\dots , x_n=t^{d_n}, \quad t\in [-1,1].\ee

\bt
If $d_1\geq 1, \ d_2 > d d_1, \dots ,  \ d_n > d d_{n-1}$, then $S$ is
${\mathcal P}_d({\mathbb R}^n)$-norming, and $N_{{\mathcal P}_d({\mathbb R}^n)}(S)\leq 2^{dd_n}\cdot\left(^d_n\right)$.
\et
\pr
Let $P(x_1,\dots, x_n)=\sum_{\vert \alpha \vert \leq d}a_\alpha x^\alpha$ be a real polynomial of degree $d$, $\alpha=(\alpha_1, \dots, \alpha_n)$ - multi-indices,
and $\vert \alpha \vert = \alpha_1+\dots+\alpha_n.$
On the curve $S$ we have $x^\alpha=t^{\beta(\alpha)},$ where
$\beta(\alpha)=\alpha_1 d_1+\alpha_2 d_2 +\dots + \alpha_n d_n.$

\bl For $\alpha' \ne \alpha''$ we have $\beta(\alpha') \ne \beta(\alpha'')$.
\el\pr Let $j\leq n$ be the largest index for which $\alpha'_j\ne \alpha''_j$,
say, $\alpha'_j < \alpha''_j$. Then
$\beta(\alpha')\leq \beta(\alpha'')-d_j+dd_{j-1}<\beta(\alpha'')$,
since $d_i$ increase, $\vert \alpha \vert \leq d$, and $d_j > d d_{j-1}$.
$\square$

\medskip

So the monomials of $P$ remain ``separated" in the univariate polynomial
$G(t)=P(\Psi(t))$ of degree $\beta\bigl((0,\dots,0,d)\bigr)=dd_n$, and their
coefficients $a_\alpha$ remain the same.  If $P$ is bounded by $1$ on $S$,
then $G$ is bounded by $1$ on $[-1,1]$. We conclude via the Chebyshev inequality
that all the coefficients of $G$ do not exceed $2^{dd_n},$ and hence the same is true for $P$. 
Finally, this implies that $P$ is bounded by $2^{dd_n}\cdot\left(^d_n\right)$
on $Q^n_1$.\qquad $\square$

\medskip

Applying to the univariate polynomial $G$ the classical Remez inequality,
or its discrete version given by Theorem 2.1 above, we immediately conclude that the corresponding subsets of the curve $S$ are  ${\mathcal P}_d({\mathbb R}^n)$-norming in ${\mathbb R}^n$.

\medskip

More general class of real analytic curves $S$ for which an analog of Theorem 2.2 is valid for certain spaces $V$ of real analytic functions with explicit bounds of norming constant
$N_{V}(S)$ is presented in \cite[Th.~2.3]{Bru2}. The role of degree of a polynomial there plays ``valency'' of an analytic function.

\subsection{Transcendental Surfaces}

Each piece of an analytic curve
\[
\Gamma:=\{\bigl(x,\phi(x)\bigr)\in {\mathbb R}^2\, :\, x\in [a,b]\subset [-1,1],\quad \sup_{[a,b]}|\phi|<1\}\subset Q^2_1
\]
is ${\mathcal P}_{d}({\mathbb R}^2)$-norming for all $d$ if $\phi$ is a transcendental function. For instance, this is
true for the curve $y=e^x$, $x\in [-1,0]$. However, it may be a delicate problem to bound explicitly the norming constants for such curves. In this section we formulate one of the results in this direction. Let us recall that an entire function $f$ on ${\mathbb C}^n$ is of order $\rho\ge 0$ if
\[
\rho=\limsup_{r\rightarrow\infty}\frac{\ln\, m_f(r)}{\ln r},\qquad {\rm where}\qquad m_f(r)=\ln\left(\sup_{\|z\|_{2}\le r} |f(z)|\right).
\]
If $\rho < \infty$, then $f$ is called of finite order.

Suppose that $f$ is a nonpolynomial entire function on ${\mathbb C}^n$ real on ${\mathbb R}^n$ and such that the hypersurface $\{(x_1,\dots, x_{n+1})\in {\mathbb R}^{n+1}\, :\, x_{n+1}=f(x_1,\dots, x_n)\}$ intersects $Q^n_1$ by a subset $Z$ of real dimension $n$. The next result follows from \cite[Th.~2.5]{Bru2} (see also \cite[Th.~1.1]{CP} for the particular case $n=1$ and $f$ being of finite positive order).
\bt\label{te3.1}
There exist a convergent to $\infty$ sequence of natural numbers $\{d_j\}$  (depending on $f$ only) and a convergent to $0$ sequence of positive numbers $\{\epsilon_j\}$ (depending on $f$, $n$ and $n$-dimensional Hausdorff measure of $Z$ only)
such that
\[
\ln\left(N_{{\mathcal P}_{d_j}({\mathbb R}^{n+1})}(Z)\right)\le d_j^{2+\epsilon_j}\quad {\rm for\ all}\quad d_j.
\]
\et
For $n=1$ this inequality is sharp in the sense that $2$ in the exponent on the right-hand side cannot be replaced by a smaller number.

In \cite[Th.~2.8]{Bru2} some sufficient conditions for $f$ are formulated under which the above inequality is valid for all polynomial degrees $d$. In particular, this is true if
\[
f=\sum_{j=1}^m p_j\cdot e^{q_j}, \quad {\rm where\ all}\quad p_j, q_j\quad {\rm are\ real\ polynomials\ on\ } {\mathbb R}^n.
\]

Also, \cite[Th.~2.5]{Bru2} deals with some other spaces $V$ of analytic functions and gives upper bounds for their norming constants $N_{V}(Z)$.\smallskip

\subsection{Topological Conditions}\label{Top.Cond}

Algebraic Geometry provides a wide variety of specific topological properties
of real algebraic sets in ${\mathbb R}^n$ of a given degree $d$. Some of them are shared also by all {\it subsets} of these algebraic sets. So if a compact set $Z\subset Q^n_1$ violates one of such properties, it cannot be contained in an algebraic set of degree $d$, and hence $Z$ is ${\mathcal P}_d({\mathbb R}^n)$-norming. Below, we illustrate this by some simple examples.

\medskip

We consider compact oriented hypersurfaces $S$ in ${\mathbb R}^n$, bounding the corresponding domains $D$. A sequence $\Sigma =\{S_1,\dots,S_m\}$ of such hypersurfaces is called ``nested" if $D_1\subset D_2 \subset \dots \subset D_m$. We define $\delta(\Sigma)$ as the minimum of the $\ell^\infty$ distances between the subsequent hypersurfaces in $\Sigma$.

\bp
Let $\Sigma = \{S_1,\dots,S_{d+1}\}$ be a nested sequence of hypersurfaces in the unit cube $Q^n_1\subset {\mathbb R}^n$. Then $S=\cup^{d+1}_{j=1} S_j$ is ${\mathcal P}_{2d}({\mathbb R}^n)$-norming, and $N_{{\mathcal P}_{2d}({\mathbb R}^n)}(S) \leq R_{2d}(\delta)=T_{2d} \left({{2-\delta}\over {\delta}}\right),$ where $\delta=\delta(\Sigma)$.
\ep
\pr
Let $P\in {\mathcal P}_{2d}({\mathbb R}^n)$ be a polynomial of degree at most $2d$ bounded by $1$
on $S$. Fix a point $x_0\in D_1$ and consider the straight line $l\subset {\mathbb R}^n$ passing through $x_0$ and
a point $y\in Q^n_1$ such that $|P(y)|=\max_{Q^n_1}|P|$. We set $l\setminus\{x_0\}:= l_1\sqcup l_2$, where each $l_i$ is an open ray with endpoint $x_0$. By assumptions, each $l_i$ crosses $S$ at not less than $d+1$ points, and  we can fix exactly one point $x_{ij}\in l_i\cap S_j$, $j=1,\dots , d+1$, so that 
\[
\|x_{i1}-x_0\|_\infty <\|x_{i2}-x_0\|_\infty <\cdots <\|x_{id+1}-x_0\|_\infty.
\]
Also, the points
$x_{i1}, x_{i2}, \dots, x_{id+1}$ on $l_i$ are separated by the $\ell^\infty$ distance
at least $\delta=\delta(\Sigma)$ from one another. This implies that the points $x_{11}$, $x_{ij}$, $2\le j\le d+1$, $i=1,2$, are separated by the $\ell^\infty$ distance
at least $\delta=\delta(\Sigma)$ from one another as well. Applying Corollary 2.2 to the set $Z$ consisting of these points and
to the restriction of $P$ to the interval $l \cap Q^n_1$ (which has the $\ell^\infty$ length
at most $2$) we get the required bound. \quad $\square$

\medskip

Let us give an example of more subtle topological restrictions.

\bt
Each compact set $Z\subset Q^2_1 \subset {\mathbb R}^2$ containing $11$ ovals out from one another is ${\mathcal P}_6({\mathbb R}^2)$-norming.
\et
\pr
If $Z$ were contained in an algebraic curve $Y$ of degree
$6$, then each oval of $Z$ would be an oval of $Y$. But by the solution of the first part of Hilbert's 16-th problem, $Y$ cannot contain $11$ ovals out from 
one another. Therefore $Z$ cannot be contained in an algebraic curve of degree $6$. \qquad $\square$

\medskip

\noindent{\bf Remark} It is not easy to give an explicit bound on $N_{{\mathcal P}_{6}({\mathbb R}^2)}(Z)$ in geometric terms. Producing such bounds (in this and similar
situations) is an important open problem. Its solution would clarify the interconnection
of topological and analytical properties of polynomials. It would also clarify some ``rigidity" properties of smooth functions appearing in the framework of
the approach developed in \cite{Yom6} (see also references therein) and aimed to transfer to several variables the classical Rolle lemma and some of its
important consequences. In particular, combining the Taylor
reminder formula with a bound for the norming constant in Theorem 2.4 would lead to the statement that any $C^7$ function $f$ on the plane,
vanishing on $Z$ as above and satisfying $\max_{x\in Q^2_1} \vert f(x) \vert=1$, must have its $7$-th derivative larger than a certain explicit constant $c_Z$ 
(depending on $Z$ only).

This last statement is directly related to the Whitney extension problem for smooth functions, especially in the form considered recently by Ch. Fefferman 
(see \cite{Fef} and references therein). We expect that Remez type inequalities can improve our understanding of the geometry of the Whitney extensions, and
plan to present some results in this direction separately.

\section{Remez-type Inequalities for Fewnomials} \label{Few} 
\setcounter{equation}{0}

\subsection{Tur\'an-Nazarov Inequality}
This is an important {\it nonlinear} version of the Remez inequality having numerous applications in Analysis, Random Functions theory, Sampling theory etc. (see
\cite{T,N,Naz.Nis.Sod,Bat.Sar.Yom} and references therein), formulated as follows:

\bt [\cite{N}]\label{turan}
Let $p(t) = \sum_{k=0}^m c_k e^{\lambda_k t}$, $t\in {\mathbb R}$, where all $c_k ,\lambda_k \in \C$, be an exponential polynomial. Let $I \subset {\mathbb R}$ be
an interval and $Z \subset I$ be a measurable subset. Then
\be \label{TN}
\sup_{I} |p| \le e^{\mu_1(I) \cdot \max_{0\le k\le m}\, |{\rm Re}\, \lambda_k |} \cdot \left( \frac{c \mu_1(I)}{\mu_1(Z) }\right)^{m} \cdot \sup_{Z} |p|,
\ee
where $\mu_1$ is the Lebesgue measure on ${\mathbb R}$ and $c>0$ is an absolute constant.
\et

This result was established, first, by Tur\'an \cite{T} for all $\lambda_k$ being pure imaginary Gaussian integers and $Z\subset I$ being an interval, and in the general
form by Nazarov \cite{N}. Its discrete version, in the spirit of \cite{Yom5}, was obtained in \cite{Fri.Yom}.\smallskip

Let us prove the multidimensional version of Theorem \ref{turan}.

\bt\label{turan1}
Let $A\subset {\mathbb R}^n$ be a $d$-dimensional affine subspace, 
$B\subset A$ be a convex body in $A$ and $Z\subset B$ be a Borel subset. Let $p(x)=\sum_{k=0}^m c_k e^{f_k(x)}$, $x\in {\mathbb R}^n$, where
all $c_k\in {\mathbb C}$ and all $f_k$ are complex-valued linear functionals on ${\mathbb R}^n$, be an exponential polynomial. Then
\be\label{TN1}
\sup_{B}|p|\le e^{\max_{0\le k\le m}\left\{\sup_{x,y\in B}{\rm Re} f_k(x-y)\right\}}\cdot\left(\frac{cd {\mathcal H}_d(B)}{{\mathcal H}_d(Z)}\right)^m\cdot\sup_Z |p|,
\ee
where ${\mathcal H}_d$ is the Hausdorff $d$-measure on ${\mathbb R}^n$.
\et
\pr
Without loss of generality we assume that $Z$ is closed and ${\mathcal H}_d(Z)>0$. Let $x_0\in B$ be such that $\sup_B |p|=|p(x_0)|$. Due to \cite[Lm.~3]{Bru.Gan}, there exists a ray $l:=\{x_0+t\cdot e\in {\mathbb R}^n\, :\, \|e\|_2=1, \ t\in {\mathbb R}_+\}$ with the endpoint $x_0$ such that
\be\label{ray}
\frac{{\mathcal H}_1(B\cap l)}{{\mathcal H}_1(Z\cap l)}\le\frac{d {\mathcal H}_d(B)}{{\mathcal H}_d(Z)} .
\ee
The restriction of $p$ to $l\cap B$ has a form
\[
\sum_{k=0}^m \left(c_k\cdot e^{f_k(x_0)}\right)e^{t\cdot f_k(e)},\qquad 0\le t\le t_0,
\]
where $t_0>0$ is such that $x_0+t_0\cdot e$ belongs to the boundary of $B$ in $A$.

Thus according to inequality (\ref{TN}),
\[
\sup_{B}|p|=\sup_{B\cap l}|p|\le e^{t_0\cdot\max_{0\le k\le m}\,|{\rm Re}\, f_k(e)|}\cdot\left(\frac{c{\mathcal H}_1(B\cap l)}{{\mathcal H}_1(Z\cap l)}\right)^m\cdot\sup_{Z\cap l} |p|.
\]
It remains to use (\ref{ray}) and note that $t_0\cdot |{\rm Re}\, f_k(e)|\le \sup_{x,y\in B}{\rm Re}\, f_k(x-y)$ and $\sup_{Z\cap l} |p|\le\sup_Z |p|$.\qquad $\Box$

\subsection{Fewnomial Remez-type inequality}

We deduce from (\ref{TN1}) the following ``fewnomial'' Remez-type inequality.

Let $({\mathbb R}^*_{+})^n\subset {\mathbb R}^n$ be the set of points with positive coordinates. If $x_1,\dots, x_n$ are coordinates on ${\mathbb R}^n$ we introduce a Riemannian metric on $({\mathbb R}^*_{+})^n$ by the formula
\[
ds^2=\frac{dx_1^2}{x_1^2}+\cdots +\frac{dx_n^2}{x_n^2}.
\]

The map $e_n: {\mathbb R}^n\rightarrow ({\mathbb R}^*_{+})^n$, $e_n((x_1,\dots, x_n))=(e^{x_1},\dots, e^{x_n})$, determines an isometry between ${\mathbb R}^n$ equipped with the Euclidean metric and $({\mathbb R}^*_{+})^n$ equipped with the above introduced Riemannian metric. Thus, the latter is a geodesically complete Riemannian manifold and geodesics there are images by $e_n$ of straight lines in ${\mathbb R}^n$. In particular, a geodesic segment joining points
$x=(x_1,\dots, x_n)$ and $y=(y_1,\dots, y_n)$ in $({\mathbb R}^*_{+})^n$ has a form
\[
x\hat{\,}^{\, t}\circ y\hat{\,}^{\, 1-t}:=(x_1^t\cdot y_1^{1-t},\dots, x_n^t\cdot y_n^{1-t}),\quad 0\le t\le 1.
\]

One easily shows that such a segment is the usual convex interval joining $x$ and $y$ if and only if there exist a partition of the set $\{1,\dots, n\}$ into
disjoint subsets $I,J$ (one of which may be $\emptyset$) and a positive number $\lambda$ such that $x_i=y_i$ for all $i\in I$ and $x_j=\lambda y_j$ for all $j\in J$.\smallskip

A subset $S\subset ({\mathbb R}^*_{+})^n$ is called {\em logarithmically convex} if for each pair of points in $S$ the geodesic segment joining them belongs to $S$. In other words, $S$ is logarithmically convex if it is the image under $e_n$ of a convex subset of ${\mathbb R}^n$. 

In general, logarithmically convex sets are not convex; the class of convex and logarithmically convex sets is relatively small (for instance, it contains $d$-dimensional rectangles in $({\mathbb R}^*_{+})^n$, $0\le d\le n$, with edges parallel to coordinate axes).

A submanifold $M\subset ({\mathbb R}^*_{+})^n$ is called {\em affine} if it is the image under $e_n$ of an affine subspace of ${\mathbb R}^n$.

For a compact subset $S\subset  ({\mathbb R}^*_{+})^n$  and a natural number $d\le n$ we define
\[
K_d(S):=\frac{\max_{\{i_1,\dots,i_d\}\in I_{d}}\left\{\max_{S}\, x_{i_1}\cdots x_{i_d}\right\}}{\min_{\{i_1,\dots ,i_d\}\in I_{d}}\left\{\min_{S}\, x_{i_1}\cdots x_{i_d}\right\}},
\]
where $I_{d}$ is the family of all $d$-point subsets of the set $\{1,\dots, n\}$. 

Clearly, $K_d(S)=K_d(t\cdot S)$ for each $t>0$, and if $L$ is a diagonal matrix with positive entries, then, $K_n(L(S))=K_n(S)$. Also,
\be\label{d1}
K_d(S)\le K_1(S)^d\quad {\rm for\ all}\quad d.
\ee

If $\alpha=(\alpha_1,\dots,\alpha_n)\in {\mathbb R}^n$ and $x=(x_1,\dots, x_n),\ y=(y_1,\dots, y_n)\in ({\mathbb R}^*_{+})^n$, we set $x^{\alpha}:=x_1^{\alpha_1}\cdots x_n^{\alpha_n}$, $|\alpha|=\alpha_1+\cdots+\alpha_n$, and $\frac{x}{y}:=\left(\frac{x_1}{y_1},\dots,\frac{x_n}{y_n}\right)$. 

\bt \label{Few.Rem}
Let $M\subset ({\mathbb R}^*_{+})^n$ be a $d$-dimensional affine submanifold, $B\subset M$ be a $d$-dimensional logarithmically convex compact subset and $Z\subset B$ be a Borel subset.
Let $p(x) = \sum_{k=0}^m c_k x^{\alpha_k}$, $x\in ({\mathbb R}^*_{+})^n$, where all $\alpha_k\in {\mathbb R}^n$ and all $c_k\in {\mathbb R}$. Then
\be \label{Few.Rem.Eq}
\sup_B |p|\le\max_{0\le k\le m}\left\{\sup_{x,y\in B}\left(\frac{x}{y}\right)^{\alpha_k}\right\}\cdot \left(\frac{cd K_d(B)\cdot{\mathcal H}_d(B)}{{\mathcal H}_d (Z)}\right)^m\cdot\sup_Z |p|.
\ee
\et

\pr
The substitution $x=e_n(u), \ u\in {\mathbb R}^n$, reduces (\ref{Few.Rem.Eq}) to a particular case of (\ref{TN1}) for the exponential polynomial $p\circ e_n$. To estimate the ratio of measures $\frac{{\mathcal H}_d(e_n^{-1}(B))}{{\mathcal H}_d (e_n^{-1}(Z))}$ in the obtained inequality,
we express the Hausdorff $d$-measure of a subset $S\subset M$ as 
\[
{\mathcal H}_d(S)=\int_{e_n^{-1}(S)}J_d(e_n|_{A})(x)d {\mathcal H}_d(x),
\]
where $A=e_n^{-1}(M)\subset {\mathbb R}^n$ is a $d$-dimensional affine subspace and $J_d(e_n|_{A})$ is
the $d$-Jacobian of the map $e_n|_{A}$, see \cite[Th.~3.2.3]{Fe}. Under a suitable affine parameterization $\phi=(\phi_1,\dots,\phi_d): {\mathbb R}^d\rightarrow A$ of $A$
we obtain
\[
{\mathcal H}_d(S)=\int_{(\phi^{-1}\circ e_n^{-1})(S)}\left(\sum_{\{i_1,\dots, i_d\}\in I_d} a_{i_1,\dots, i_d}^2\cdot e^{2\phi_{i_1}(x)+\cdots+2\phi_{i_d}(x)} \right)^{\frac 12}d\mu_d(x),
\]
where $a_{i_1,\dots, i_d}$ is the determinant of the $d\times d$ matrix defined by the $i_1,\dots, i_d$ rows of the linear part of $\phi$. 
The latter implies for $S':=(e_n\circ \phi)^{-1}(S)$
\[
\begin{array}{l}
\displaystyle \left(\sum_{\{i_1,\dots, i_d\}\in I_d} a_{i_1,\dots, i_d}^2\right)^{\frac 12}\cdot \min_{\{i_1,\dots,i_d\}\in I_{d}}\left\{\min_{S}\, x_{i_1}\cdots x_{i_d}\right\}\cdot \mu_d(S')\le
{\mathcal H}_d(S)
\\
\displaystyle
\le \left(\sum_{\{i_1,\dots, i_d\}\in I_d} a_{i_1,\dots, i_d}^2\right)^{\frac 12}\cdot \max_{\{i_1,\dots,i_d\}\in I_{d}}\left\{\max_{S}\, x_{i_1}\cdots x_{i_d}\right\}\cdot \mu_d(S').
\end{array}
\]
In turn, this yields
\[
\frac{{\mathcal H}_d(e_n^{-1}(B))}{{\mathcal H}_d (e_n^{-1}(Z))}=\frac{{\mu}_d((e_n\circ\phi)^{-1}(B))}{{\mu}_d ((e_n\circ\phi)^{-1}(Z))}\le\frac{K_d(B)\cdot{\mathcal H}_d(B)}{{\mathcal H}_d (Z)}.\qquad\Box
\]
\smallskip

If $B=\{x=(x_1,\dots, x_n)\in ({\mathbb R}^*_{+})^n\, :\, a_i< x_i< b_i,\ 1\le i\le n\}$ is an $n$-dimensional rectangle, then since
${\mathcal H}_n(e_n^{-1}(B))=\ln\left(\frac{b_1}{a_1}\right)\cdots \ln\left(\frac{b_n}{a_n}\right)$,
one obtains easily from the proof that the constant in (\ref{Few.Rem.Eq}) can be replaced by a smaller one, that is for $Z\subset B$ and $p$ as above, and $a:=(a_1,\dots, a_n)$, $b:=(b_1,\dots, b_n)$,
\be\label{rectangle}
\sup_B |p|\le\left(\max_{0\le k\le m}\left(\frac{b}{a}\right)^{\alpha_k}\right)\cdot \left(\frac{cn\cdot\prod_{i=1}^n \left(b_i\ln\left(\frac{b_i}{a_i}\right)\right)}{\mu_n(Z)}\right)^m\cdot\sup_Z |p|.
\ee

Also, using inequality (\ref{d1}) we obtain:
\bc
Under the assumptions and in notations of Theorem \ref{Few.Rem},
\be\label{e3.6}
\sup_B |p|\le  K_1(B)^{m+\max_{0\le k\le m}|\alpha_k|}\cdot\left(\frac{cd \cdot{\mathcal H}_d(B)}{{\mathcal H}_d (Z)}\right)^m\cdot\sup_Z |p|.
\ee
\ec
\smallskip 

An important feature of inequality (\ref{e3.6}) is that while the degree ${\rm deg}\, p:=\max_{0\le k\le m}|\alpha_k|$ of $p$ enters the constant of the inequality
as the exponent of a certain geometric characteristic of $B$, the exponent of $\frac{{\mathcal H}_d(B)}{{\mathcal H}_d (Z)}$ is $m$, i.e. it depends only on the number 
of terms of $p$, but not on its degree. 

\medskip

If $B\subset  ({\mathbb R}^*_{+})^n$ in Theorem \ref{Few.Rem} is of the form $x_0+Q^n_1$, $x_0\in ({\mathbb R}^*_{+})^n$, then one easily deduces from inequality (\ref{Few.Rem.Eq})
its discrete version by replacing the Hausdorff $n$-measure of $Z$ in (\ref{Few.Rem.Eq})
by a fewnomial version of the metric span of $Z$. It is defined as in Definition \ref{span} above, with
the polynomial $M_{n,d}$ replaced by its ``fewnomial'' version (see \cite{Fri.Yom1}). Since the expression is rather cumbersome, we do not state this result explicitly in full generality, restricting ourselves to the particular case of univariate polynomials. 

\bt \label{Rem.Few.Discr}
Let $[a,b]\Subset \mathbb R^*_+$ and $Z\subset [a,b]$ be a measurable subset. Suppose $p(x)=\sum_{k=0}^m c_k x^{n_{k}}$, $x\in {\mathbb R}$, where $0\le n_0<\cdots< n_m$ are nonnegative integers, is a real polynomial of degree at most $n_m$. Then
\[
\sup_{[a,b]} |p|\le \left(\frac{b}{a}\right)^{n_m}\cdot\left(\frac{c\left(b\ln\left(\frac{b}{a}\right)\right)}{\omega_{m,1}(Z)}\right)^m\sup_Z |p|,
\]
where $\omega_{m,1}(Z)$ is the metric span of $Z$.
\et
\pr
The proof repeats literally the proof of Theorem 2.1 above presented in \cite{Yom5}. The only property required in the proof is that, according to
the Descartes rule, for every $C\in {\mathbb R}$, the polynomial $p-C$ has at most $m$ positive roots. \qquad $\square$

\smallskip

As a simple example of sets $Z$ satisfying fewnomial Remez-type inequalities, we consider a nested sequence of hypersurfaces $\Sigma = \{S_1,\dots,S_{m+1}\}$, as in
Section \ref{Top.Cond} above. Assume that $\Sigma \subset Q^n_R \setminus Q^n_\rho\subset Q^1_n,$ where $Q^n_s\subset {\mathbb R}^n$ stands for the closed $\ell^\infty$ ball (-``cube'') of radius $s$ centered at
$0 \in {\mathbb R}^n$. As above, let $\delta=\delta(\Sigma)$ denote the minimal $\ell^\infty$ distance between $S_j$.

Consider the family of monomials ${\mathcal F}=\{x^{\alpha_k}\}_{0\le k\le m},$ where all $\alpha_k \in {\mathbb Z}_{+}^n$ and
$N:=\max_{0\le k\le m}\,|\alpha_k|$. The linear space $V_{\mathcal F}$ generated by ${\mathcal F}$ consists of 
multivariate polynomials of the form
\be \label{feww}
P(x)=\sum_{k=0}^m c_k x^{\alpha_k}, \ c_k \in {\mathbb R}, \ x \in {\mathbb R}^n.
\ee

\bt \label{multi.few}
The set $S=\cup^{m+1}_{j=1} S_j$ is an $V_{\mathcal F}$-norming and for each
$P\in V_{\mathcal F}$, 
\be \label{multi.few.eq}
\sup_{Q^n_R} |P| \le \left({R\over \rho}\right)^N \cdot \left( \frac{c R \ln (R/\rho)}{\delta}\right)^{m} \cdot \sup_{S} |P|.
\ee
\et
\pr
Consider a ray $l= \{tv, \ v \in {\mathbb R}^n,\ ||v||_\infty =1, \ t\ge 0\}\subset {\mathbb R}^n$ with endpoint $0$ passing through
a point $y\in Q^n_R$ such that $|P(y)|=\max_{Q^n_R}|P|$, and apply Theorem
\ref{Rem.Few.Discr} to the interval $B_l= \{tv\in l\, :\, \rho \leq t \leq R\}$ and the subset $Z_l=S \cap l$.
By assumptions, $B_l$ crosses $S$ at not less than $m+1$ points, and we can fix exactly one point $x_j\in l\cap S_j$, $j=1,\dots , m+1$, so that
the points $x_{1}, x_{2}, \dots, x_{m+1}$ on $l$ are ordered and separated by the $\ell^\infty$ distance
at least $\delta=\delta(\Sigma)$ from one another. Applying the estimate of the metric span in Corollary
\ref{Discr.Rem1} above, we conclude that $\omega_{m,1}(Z_l)\geq \delta$. This gives us the required bound. \qquad $\square$

\section{Lipschitz Continuity of the Norming Constant}\label{Sec3}
\setcounter{equation}{0}
Let $(X,d)$ be a metric space. A real function $f$ on $X$ is said to belong to the space ${\rm Lip}(X)$ if 
\[
L_f:=\sup_{x\ne y}\frac{|f(x)-f(y)|}{d(x,y)}<\infty .
\] 
In this case, $L_f$ is called the {\em Lipschitz constant of} $f$.
\smallskip

Let $\omega$ be an increasing concave function on ${\mathbb R}_+$, equal to $0$ at $0$, and such that $\lim_{t\rightarrow\infty}\omega(t)=\infty$. One easily checks that $d_\omega(x,y):=\omega(\|x-y\|_\infty)$, $x,y\in {\mathbb R}^n$, is a metric on 
${\mathbb R}^n$ compatible with the standard topology. Note that $C(Q^n_1)=\bigcup_{\omega}{\rm Lip}_{d_\omega} (Q^n_1)$, where the union is taken over all possible $\omega$.\smallskip

Let $V\subset {\rm Lip}_{d_\omega} (Q^n_1)$, for some $\omega$ as above, be a finite-dimensional space.
\bp\label{Markov}
There exists a constant $M>0$ such that for each $f\in V$ 
\be\label{eqmarkov}
L_f\le M\sup_{Q^n_1}|f|.
\ee
\ep
\pr Let ${\mathcal F}=\{f_1,\dots, f_l\}$, $l:={\rm dim}\, V$,  be a basis in $V$. For $f=\sum_{i=1}^l a_i f_i\in V$, we set $\|f\|_1:=\sum_{i=1}^l |a_i|$. Since $V$ is $l$-dimensional, $\|\cdot\|_1$ is a norm
on $V$ equivalent to the norm induced from $C(Q^n_1)$. In particular, for some $\tilde c>0$ and all $x,y\in Q^n_1$, $x\ne y$, we have
\[
\frac{|f(x)-f(y)|}{d_\omega(x,y)}\le \left(\max_{1\le i\le l} \,L_{f_i}\right)\cdot\|f\|_1\le  \left(\max_{1\le i\le l} \,L_{f_i}\right)\cdot\tilde c\cdot\sup_{Q^n_1}|f|.\qquad \Box
\]
\smallskip

The optimal constant $M_V$ in (\ref{eqmarkov}) is called the {\em Markov constant of}
$V$.  For instance, the classical A.~Markov polynomial inequality implies that $M_{V}=d^2 n$ if $V={\mathcal P}_d({\mathbb R}^n)$ and $\omega(t):=t$. In turn, the classical Bernstein inequality implies that $M_{V}=\pi d n$ if  $V$ is the space of trigonometric polynomials of degree $d$ on ${\mathbb R}^n$ of period $2$ in each coordinate and $\omega(t)=t$.
The constant $M_{V}$ can be effectively estimated applying the Gram--Schmidt process to a basis ${\mathcal F}=\{f_1,\dots, f_l\}$ in $V$ considered in a suitable space $L^2(\mu)$ on $Q^n_1$:
\bp
Suppose all $f_i\in {\rm Lip}_{d_\omega} (Q^n_1)$ and the family $\mathcal F$ is orthonormal with respect to a regular Borel measure $\mu$ on $Q^n_1$. Then
\[
M_{V}\le \left(\max_{1\le i\le l} \,L_{f_i}\right)\cdot \sqrt{l}\cdot\sqrt{\mu(Q_n^1)}.
\]
\ep

\pr For $f=\sum_{i=1}^l a_i f_i\in V$ we have
\[
\|f\|_1\leq\sqrt{l}\cdot\left(\sum_{i=1}^l |a_i|^2\right)^{\frac 12}=\sqrt{l}\cdot\left(\int_{Q^n_1}f^2 d\mu\right)^{\frac 12}\le\sqrt{l}\cdot\sqrt{\mu(Q_n^1)}\cdot \sup_{Q^n_1}|f|.
\]

This and the argument of the proof of Proposition \ref{Markov} give the required inequality.\qquad $\Box$\smallskip

Let ${\mathcal K}_n$ be the set of all closed subsets of $Q^n_1$ equipped with the Hausdorff metric $d_H$: if $K_0,K_1\in {\mathcal K}_n$, then
\[
d_H(K_0,K_1):=\max_{i=0,1}\left\{\sup_{y\in K_i} \inf_{x\in K_{1-i}} \|x-y\|_\infty\right\}.
\]

It is well known that $({\mathcal K}_n, d_H)$ is a compact metric space. Let us consider ${\mathcal K}_n$ with the metric $d_{\omega H}:=\omega\circ d_H$. Then one checks easily that $({\mathcal K}_n, d_{\omega H})$ is compact as well and that the metrics $d_H$ and $d_{\omega H}$ determine the same topology on ${\mathcal K}_n$.

Our main result is the following Lipschitz continuity property of norming constants $N_{V}(Z)$, $Z\in {\mathcal K}_n$.

\bt
The function $\frac{1}{N_{V}}\in {\rm Lip}_{d_{\omega H}}({\mathcal K}_n)$ and its Lipschitz constant $L_{\frac{1}{N_{V}}}\le M_{V}$.
(Here we define $\frac{1}{N_{V}(Z)}=0$ for $Z$ not $V$-norming.)
\et
\pr Let $Z_1,Z_2\in {\mathcal K}_n$ be $V$-norming sets.  Assume without loss of generality that $N_{V}(Z_1)\ge N_{V}(Z_2)$. Suppose that $f\in V$ is such that $\sup_{Z_1}|f|=1$ and $\sup_{Q^n_1}|f|=N_{V}(Z_1)$. For each $z_2\in Z_2$ we choose $z_1\in Z_1$ so that
$\|z_2-z_1\|_\infty\le d_{H}(Z_2,Z_1)$. Then we have
\[
|f(z_2)|\le |f(z_2)-f(z_1)|+|f(z_1)|\le M_{V}\cdot d_{\omega H}(Z_1,Z_2)\cdot N_{V}(Z_1)+1.
\]
Thus, by the definition of $N_{V}(Z_2)$,
\[
N_{V}(Z_1)=\sup_{Q^n_1}|f|\le (M_{V}\cdot d_{\omega H}(Z_1,Z_2)\cdot N_{V}(Z_1)+1)\cdot N_{V}(Z_2).
\]
This implies the required statement:
\[
N_{V}(Z_1)-N_{V}(Z_2)\le M_{V}\cdot d_{\omega H}(Z_1,Z_2)\cdot N_{V}(Z_1)\cdot N_{V}(Z_2).
\]

Further, assume that $Z_1$ is not $V$-norming while $Z_2$ is. Then due to Proposition \ref{prop1} there exists a function $f\in V$ such that $f|_{Z_1}=0$ and
$\sup_{Q^n_1}|f|=1$. Arguing as above (with $N_{V}(Z_1)$ replaced by $1$) we obtain
\[
1=\sup_{Q^n_1}|f|\le M_{V}\cdot d_{\omega H}(Z_1,Z_2)\cdot N_{V}(Z_2),
\]
that is
\[
\frac{1}{N_{V}(Z_2)}-\frac{1}{N_{V}(Z_1)}\le M_{V}\cdot d_{\omega H}(Z_1,Z_2).\qquad \Box
\]

\bc If $Z\in {\mathcal K}_n$ is $V$-norming, then each $Y$ in the open ball of radius $\frac{1}{M_{V}N_{V}(Z)}$ with center at $Z$ is $V$-norming and
\[
N_{V}(Y)\le  \frac{N_{V}(Z)}{1-M_{V}\cdot N_{V}(Z)\cdot d_{\omega H}(Z,Y)}.
\]
\ec

\pr From the previous theorem we obtain
\[
\frac{1}{N_{V}(Y)}=\frac{1}{N_{V}(Z)}-\left(\frac{1}{N_{V}(Z)}-\frac{1}{N_{V}(Y)}\right)\ge \frac{1}{N_{V}(Z)}-M_{V}\cdot d_{\omega H}(Z,Y)>0.
\]
Passing here to reciprocals we get the required statement.\qquad $\Box$

\smallskip

In particular, the set of non $V$-norming sets is a closed subset of ${\mathcal K}_n$. Also, one easily shows that if $V$ consists of analytic functions, then the set of all
$V$-norming sets is open and dense in ${\mathcal K}_n$.

\medskip
\medskip

\bibliographystyle{amsplain}

\begin{thebibliography}{10}

\bibitem{AT} H. Alexander and B. A. Taylor, Comparison of two capacities in ${\mathbb C}^n$, {\sl Math. Z.} {\bf 186} (1984), 407--417.

\bibitem{And} V. V. Andrievskii, Local Remez-type inequalities for
exponentials of a potential on a piecewise analytic arc,
{\sl Journal d'Analyse Math.}  {\bf 100}, No. 1 (2006), 323--336.

\bibitem{Ar} N. Aronszajn, A unique continuation theorem for solutions of elliptic partial differential equations
of second order, {\sl J. Math. Pures Appl.} {\bf 36} (1957), 235--249.

\bibitem{Bat.Sar.Yom} D. Batenkov, N. Sarig and Y. Yomdin, Accuracy of algebraic Fourier reconstruction for shifts of several signals,
submitted, arXiv:1311.3468.

\bibitem{Ber2} S. Bernstein, Sur la limitation des valeurs d'un
polynome $P_n(x)$ de degr\'e $n$ sur tout un segment par ses
vleurs en $n+1$ points du segment, {\sl Isvestiya AN SSSR}
(1931), 1025--1050.

\bibitem{BE} P.~B.~Borwein and T. Erd\'{e}lyi, Generalizations of M\"{u}ntz's theorem via a Remez-type inequality
for M\"{u}ntz spaces, {\sl J. Amer. Math. Soc.}
{\bf 10} (1997), 327--349.

\bibitem{BDSV}
L. Bos, S. De Marchi, A. Sommariva and M. Vianello,  Weakly admissible meshes and discrete extremal sets, {\sl Numer. Math. Theory Methods Appl.} {\bf 4} (2011), 1--12.

\bibitem{Bru} A. Brudnyi, On covering numbers of sublevel sets of analytic
functions, {\sl J.
Approx. Theory} {\bf 162}
(2010), 72--93.

\bibitem{Bru-1} A. Brudnyi, On local behavior of analytic functions, {\sl J. Funct.
Anal.} {\bf 169} (1999), no. 2, 481--493.

\bibitem{Bru0} A. Brudnyi, Local inequalities for plurisubharmonic functions, {\sl Ann. Math.} {\bf 149} (1999),
511--533.



\bibitem{Bru2} A. Brudnyi,  On local behavior of holomorphic functions along complex submanifolds of ${\mathbb C}^N$, {\sl Invent. Math.} {\bf 173} (2) (2008), 315--363.


\bibitem{Bru.Bru} A. Brudnyi and Yu. Brudnyi, Remez type inequalities
and Morrey-Campanato spaces on Ahlfors regular sets, {\sl Contemporary
Mathematics} {\bf 445} (2007), 19--44.

\bibitem{BB} A. Brudnyi and Yu. Brudnyi, Methods of Geometric Analysis in Lipschitz extension and trace
problems, Volumes I, II, {\sl Birkh\"{a}user}, 2012.

\bibitem{Bru.Gan} Yu. Brudnyi and M. Ganzburg, On an extremal problem
for polynomials of $n$ variables, {\sl Math. USSR Izv.} {\bf 37}
(1973), 344--355.

\bibitem{CP} D. Coman and E. A. Poletsky, Transcendence measures and algebraic growth of
entire functions, {\sl Invent. Math.} {\bf 170} (1) (2007), 103--145.

\bibitem{Cop.Riv} D. Coppersmith and T. J. Rivlin, The growth of
polynomials bounded at equally spaced points, {\sl SIAM J. Math.
Anal.} {\bf 23} (1992), no. 4, 970--983.


\bibitem{DR} R. Dudley and B. Randol, Implications of pointwise bounds on polynomials,
{\sl Duke Math. J.} {\bf 29} (1962), 455--458.

\bibitem{Erd} T. Erdelyi, Remez-type inequalities and their
applications, {\sl J. Comp. Appl. Math.} {\bf 47} (1993), 167--209.

\bibitem{Fav} J. Favard, Sur l'interpolation, {\sl Bull. de la S.
M. F.} {\bf 67} (1939), 103--113.

\bibitem{Fe}
H.~Federer, Geometric Measure Theory, {\sl Springer}, 1969.

\bibitem{Fef}  Ch. Fefferman, A. Israel and G. Luli, Sobolev extension by linear operators, {\sl J. Amer. Math. Soc.} {\bf 27}  
(2014),  no. 1, 69--145. 

\bibitem{Fri.Yom} O. Friedland and Y. Yomdin, An observation on Turan-Nazarov inequality, to appear, arXiv:1107.0039.

\bibitem{Fri.Yom1} O. Friedland and Y. Yomdin, Vitushkin-type theorems, to appear, arXiv:1302.5373.

\bibitem{HS} R. Hardt and L. Simon, Nodal sets for solutions of elliptic equations, {\sl J. Differential Geometry}
{\bf 30} (1989), 505--522.

\bibitem{Kli} M. Klimek, Pluripotential theory, {\sl Oxford University Press}, New York, 1991.

\bibitem{Kol.Tih} A. N. Kolmogorov and V. M. Tihomirov,
$\varepsilon$-entropy and $\varepsilon$-capacity of sets in
functional space, {\sl Amer. Math. Soc. Transl.} {\bf  17}
(1961), 277--364.

\bibitem{Kro}
A. Kro\'{o}, Bernstein type inequalities on star-like domains in ${\mathbb R}^{d}$
with application to norming sets, {\sl Bull. Math. Sci.} {\bf 3} (2013), 349--361.

\bibitem{Lab} D. Labutin, Pluripolarity of sets with small Hausdorff measure, {\sl Manuscripta math.} {\bf 102}
(2000), 163--167.

\bibitem{Lor} R. Lorentz, Multivariate Birkhoff Interpolation,
{\sl Springer-Verlag}, 1992.


\bibitem {N} F.~L. Nazarov, Local estimates for exponential polynomials and their applications to inequalities of the uncertainty principle
type, {\sl Algebra i Analiz} {\bf 5} (1993), no. 4, 3--66; translation in {\sl St. Petersburg Math. J.} {\bf 5} (1994), no. 4, 663--717.

\bibitem{Naz.Nis.Sod} F. Nazarov, A. Nishry and M. Sodin, Log-integrability of Rademacher Fourier series, with applications to random analytic functions,
arXiv:1301.0529.


\bibitem{Rah} E. A. Rakhmanov, Bounds for polynomials with a unit
discrete norm, {\sl Ann. of Math.} (2) {\bf 165} (2007), no. 1,
55--88.

\bibitem{Rem} E. J. Remez, Sur une propriete des polynomes de
Tchebycheff, {\sl Comm. Inst. Sci. Kharkov} {\bf 13} (1936), 93--95.



\bibitem{TL} B. A. Taylor and N. Levenberg, Comparison of capacities in ${\mathbb C}^n$, {\sl Lecture Notes in Math.}
{\bf 1094} (1984), 162--172.

\bibitem{Tim} A. F. Timan, Theory of approximation of functions of real variable, {\sl Pergamon
Press}, Oxford, 1963.

\bibitem {T} P. Tur\'an,  Eine neue Methode in der Analysis und deren Anwendungen, {\sl Akad\'emiai Kiad\'o}, Budapest, 1953. 

\bibitem{Vit1} A. G. Vitushkin, O mnogomernyh variaciyah, {\sl Gosudarstv. Izdat. Tehn.-Teor. Lit.}, Moscow, 1955 (Russian).

\bibitem{W} K. Weierstra\ss, {\sl Zu Hrn. Lindemann's Abhandlung:``\"{U}ber die Ludolph'sche Zahl''}, Sitzungber. K\"{o}nigl.
Preuss. Akad. Wissensch. zu Berlin {\bf 2} (1885), 1067--1086.





\bibitem{Yom5} Y. Yomdin, Remez-type inequality for discrete sets, {\sl Isr. J. Math.} {\bf 186} (2011), 45--60.

\bibitem{Yom6} Y. Yomdin, Generalized Remez inequality for $(s,p)$-valent functions, arXiv:1102.2580.

\bibitem{Yom7} Y. Yomdin, Remez-type inequality for smooth functions, to appear, arXiv:1306.3641.


\bibitem{Zer1} A. Zeriahi, A minimum principle for plurisubharmonic
functions, {\sl Indiana Univ. Math. J.} {\bf 56} No. 6 (2007),
2671--2696.


\end{thebibliography}

\end{document}